\documentclass[12pt,a4paper]{amsart}
\usepackage[top=1.15in, bottom=1.15in, left=1.17in, right=1.17in]{geometry}
\usepackage{amssymb,amsmath}
\usepackage{pdflscape}
\usepackage{amscd}
\usepackage{fancyhdr}
\usepackage{tikz}
\usetikzlibrary{matrix,arrows,cd}
\usepackage{verbatim}
\usepackage{todonotes}

\DeclareMathOperator{\Hom}{Hom}
\DeclareMathOperator{\Ext}{Ext}
\DeclareMathOperator{\End}{End}

\newtheorem{theorem}{Theorem}[section]
\newtheorem{corollary}[theorem]{Corollary}
\newtheorem{lemma}[theorem]{Lemma}
\newtheorem{proposition}[theorem]{Proposition}

\title{Fano quiver moduli}
\author{Hans Franzen \and Markus Reineke \and Silvia Sabatini}
\address{
Hans Franzen\\ Ruhr-University Bochum\\ Faculty of Mathematics\\
Universit\"ats\-strasse 150\\
44780 Bochum (Germany)}
\email{hans.franzen@rub.de}
\address{ 
Markus Reineke\\ Ruhr-University Bochum\\ Faculty of Mathematics\\
Universit\"ats\-strasse 150\\
44780 Bochum (Germany)}
\email{Markus.Reineke@ruhr-uni-bochum.de}
\address{
Silvia Sabatini\\ University of Cologne\\ Mathematical Institute\\
Weyertal 86-90\\
50931 Cologne (Germany)}
\email{sabatini@math.uni-koeln.de}

\begin{document}

\begin{abstract} We exhibit a large class of quiver moduli spaces which are Fano varieties, by studying line bundles on quiver moduli and their global sections in general, and work out several classes of examples, comprising moduli spaces of point configurations, Kronecker moduli, and toric quiver moduli.
\end{abstract}
\maketitle
\parindent0pt

\section{Introduction}

Moduli spaces of stable quiver representations, introduced in \cite{King}, form a particularly tractable class
of moduli spaces, since they are easily constructed (as a GIT quotient for a linear representation of a reductive
algebraic group), and since they arise from classification problems for linear algebra data.\\[1ex]
In the (numerically characterized) case of fine quiver moduli spaces for acyclic quivers,
which are smooth projective varieties (and rational by \cite{SchB}), several favorable geometric
properties of cohomological nature are well-established:

Their singular cohomology is algebraic \cite{KW}, they are polynomial-count
\cite{RHN}, the Betti numbers (and the counting polynomial) can be determined explicitly \cite{RHN}, and
their (torus-equivariant) Chow rings are tautologically generated and presented \cite{F1,F2}.\\[1ex]
In the present work, we add to these properties one of algebraic-geometric nature: under a rather mild
genericity assumption, we exhibit for each indivisible dimension vector for an acyclic quiver a class
of stabilities (the canonical chamber) for which the moduli space is a Fano variety, of known dimension,
Picard rank, and index (see Theorem \ref{mainthm}).\\[1ex]
A similar result was announced in \cite{Fei}, and can also be derived from the methods of \cite{Hal}. The second named author was informed by L.~Hille that he derived a similar result in \cite{Hil}. In the special case of toric quiver moduli spaces, the Fano property is established in \cite{AH}.\\[1ex]
Our result hopefully allows a fruitful interaction between the theory of Fano
varieties and the theory of quiver moduli: on the one hand, it adds classes of arbitrarily
high-dimensional, although certainly very special, Fano varieties to the existing classes of examples. In particular, it seems likely that deep properties of Fano varieties, for example the Mukai conjecture and its various generalizations, can be verified for this class (in the spirit of such a verification for toric varieties \cite{Ca06}, horospherical varieties \cite{Pas10} and symmetric varieties \cite{GH}). 
On the other hand, the particularly rigid properties of Fano varieties and the refined invariants
available for them should facilitate the study of the geometry of these quiver moduli.\\[1ex]
After reviewing the construction of quiver moduli spaces and their basic geometric properties in
Section \ref{s1}, we turn to the study of their line bundles in Section \ref{s2}: first recalling the
construction of tautological bundles, we determine (under the assumption of ample stability)
their Picard group, and identify one chamber of their ample cone. This allows us to establish the
Fano property in Section \ref{s3}, by performing a Chern class computation for the class
of the tangent bundle, which readily leads to the relevance of the canonical stability chamber and to the main theorem.\\[1ex]
The second main focus of the present paper is on illustrating the generality of our assumptions by establishing several classes of examples in Section \ref{s4}:
we first identify all fine moduli spaces of ordered point configurations in projective spaces as Fano
varieties, also discussing a subtle example related to the Segre cubic. Then we use results of
\cite{RS} to prove that all fine Kronecker moduli spaces are Fano (of non-trivial index), and verify
the Kobayashi-Ochiai theorem for them purely numerically. Moreover, we exhibit a class of Fano quiver moduli spaces of arbitrary Picard rank and index, which can also be interpreted as certain torus quotients of Grassmannians, for which we verify the Mukai conjecture (see Proposition \ref{vmukai}). We finally parametrize all Fano toric quiver moduli arising
from our result, which allows us to exhibit particular examples of two- and three-dimensional
Fano quiver moduli spaces.\\[2ex]
{\bf Acknowledgments:} The first named author would like to thank A.\ Brecan, A.\ Huckleberry, and T.\ Peternell for interesting discussions on the subject. The second named author would like to thank P.\ Belmans for helpful
explanations on Fano varieties and the Fanography \cite{Fanog}, for carefully reading a previous version, and for suggesting various improvements, as well as K. Martinez for pointing out an imprecision in the verification of the Fano property for moduli spaces of point configurations.
The authors are supported by the DFG SFB / Transregio 191 ``Symplektische Strukturen in Geometrie, Algebra und Dynamik''. The second named author is supported by the DFG GRK 2240 ``Algebro-geometrische Methoden in Algebra, Arithmetik und Topologie''.

\section{Recollections on quiver moduli}\label{s1}

\subsection{Construction of quiver moduli}

Let $Q$ be a finite quiver with set of vertices $Q_0$ and arrows written $\alpha:i\rightarrow j$. We assume $Q$ to be acyclic, that is, $Q$ has no oriented cycles. We define the Euler form $\langle\_,\_\rangle_Q$ of $Q$ on $\mathbb{Z}Q_0$ (with natural basis ${\bf i}$ for $i\in Q_0$) defined by $$\langle{\bf d},{\bf e}\rangle=\sum_{i\in Q_0}d_ie_i-\sum_{\alpha:i\rightarrow j}d_ie_j$$
for ${\bf d}=\sum_{i\in Q_0}d_i{\bf i}$ (and similarly for ${\bf e}$).\\[1ex]
Let ${\bf d}\in\mathbb{N}Q_0$ be a dimension vector. Fixing vector spaces $V_i$ of dimension $d_i$ for all $i\in Q_0$, let $$R_{\bf d}(Q)=\bigoplus_{\alpha:i\rightarrow j}{\rm Hom}_\mathbb{C}(V_i,V_j)$$ be the space of complex representations of $Q$ of dimension vector ${\bf d}$, on which the group $$G_{\bf d}=\prod_{i\in Q_0}{\rm GL}(V_i)$$ naturally acts via base change in all $V_i$. Note that the diagonally embedded scalar group $\mathbb{C}^*\subset G_{\bf d}$ acts trivially.\\[1ex]
We consider the dual group $(\mathbb{Z}Q_0)^*$ with its natural basis elements $i$ taking the $i$-th component of a dimension vector. Let $\Theta\in(\mathbb{Z}Q_0)^*$ be a linear form such that $\Theta({\bf d})=0$, and define a point $(f_\alpha)_\alpha\in R_{\bf d}(Q)$, considered as a representation $V=((V_i)_i,(f_\alpha)_\alpha)$ of $Q$, to be $\Theta$-semistable (resp.~$\Theta$-stable) if $\Theta({\bf e})\leq 0$ (resp.~$\Theta({\bf e})<0$) for all dimension vectors ${\bf e}$ of proper non-zero subrepresentations of $V$.\\[1ex]
Denoting by $R_{\bf d}^{\Theta-{\rm (s)st}}(Q)$ the Zariski open (semi-)stable locus in $R_{\bf d}(Q)$, we consider the GIT quotient
$$M_{\bf d}^{\Theta-{\rm sst}}(Q)=R_{\bf d}^{\Theta-{\rm sst}}(Q)//G_{\bf d}$$
and the geometric quotient $$M_{\bf d}^{\Theta-{\rm st}}(Q)=R_{\bf d}^{\Theta-{\rm st}}(Q)/G_{\bf d}.$$
We briefly review their construction following \cite{King}: for a character $\chi$ of $G_{\bf d}$, we consider the space $\mathbb{C}[R_{\bf d}(Q)]^{G_{\bf d},\chi}$ of $\chi$-semi-invariant polynomial functions on $R_{\bf d}(Q)$, and the graded ring
$$\mathbb{C}[R_{\bf d}(Q)]^{G_{\bf d}}_\chi=\bigoplus_{N\geq 0}\mathbb{C}[R_{\bf d}(Q)]^{G_{\bf d},N\chi}.$$
To the stability $\Theta$, we associate the character
$$\chi_\Theta((g_i)_i)=\prod_{i\in Q_0}\det(g_i)^{-\Theta({\bf i})}$$
(the sign change resulting from our definition of $\Theta$-semistability in contrast to the convention of \cite{King}). Then
$$M_{\bf d}^{\Theta-{\rm sst}}(Q)={\rm Proj}(\mathbb{C}[R_{\bf d}(Q)]^{G_{\bf d}}_{\chi_\Theta}]).$$
More explicitly, choosing $N$ big enough so that $ \mathbb{C}[R_{\bf d}(Q)]^{G_{\bf d}}_{\chi_\Theta}$ is generated by homogeneous elements $f_0,\ldots,f_s$ in $\mathbb{C}[R_{\bf d}(Q)]^{G_{\bf d},N\chi_\Theta}$, the moduli space $M_{\bf d}^{\Theta-{\rm sst}}(Q)$ equals the image of the map $$f=(f_0:\ldots:f_s):R_{\bf d}^{\Theta-{\rm sst}}(Q)\rightarrow\mathbb{P}^s.$$
Finally, the moduli space $M_{\bf d}^{\Theta-{\rm st}}(Q)$ equals the (open) image of $R_{\bf d}^{\Theta-{\rm st}}(Q)$ under the map $f$.\\[1ex]
This moduli space, if non-empty, is called the moduli space of $\Theta$-stable representations of $Q$ of dimension vector ${\bf d}$; it is a connected complex algebraic manifold of complex dimension $1-\langle{\bf d},{\bf d}\rangle$. Its points naturally correspond to isomorphism classes $[V]$ of $\Theta$-stable representations $V$ of $Q$ of dimension vector ${\bf d}$. Furthermore, the quotient map $$R_{\bf d}^{\Theta-{\rm st}}(Q)\rightarrow M_{\bf d}^{\Theta-{\rm st}}(Q)$$ is a a principal bundle for the group $PG_{\bf d}=G_{\bf d}/\mathbb{C}^*$.\\[1ex]
We call the dimension vector ${\bf d}$ indivisible if ${\rm gcd}(d_i\, :\, i\in Q_0)=1$, and $\Theta$-coprime if $\Theta({\bf e})\not=0$ for all non-zero proper ${\bf e}\leq{\bf d}$. Then $\Theta$-coprimality implies indivisibility, and, conversely, an indivisible dimension vector is $\Theta$-coprime for a sufficiently generic choice of $\Theta$. In case ${\bf d}$ is $\Theta$-coprime, the $\Theta$-stable and the $\Theta$-semistable locus in $R_{\bf d}(Q)$ coincide, thus $M_{\bf d}^{\Theta-{\rm (s)st}}(Q)$ is a connected compact complex algebraic manifold.

\subsection{Stabilities}

We denote by $\{{\bf d},{\bf e}\}=\langle{\bf d},{\bf e}\rangle-\langle{\bf e},{\bf d}\rangle$ the antisymmetrized Euler form of $Q$ and use it to define the stability $\{{\bf d},\_\}$ of $Q$, which we call the canonical stability for ${\bf d}$. We can reformulate \cite[Theorem 6.1]{Scho} as follows: there exists a $\{{\bf d},\_\}$-stable representation of dimension vector ${\bf d}$ if and only if there exists a representation $V$ of dimension vector ${\bf d}$ with trivial endomorphism ring ${\rm End}_Q(V)=\mathbb{C}$.\\[1ex]
For an indivisible dimension vector ${\bf d}$, we consider the abelian group ${\rm Stab}({\bf d})=\{\Theta\in(\mathbb{Z}Q_0)^*\, :\, \Theta({\bf d})=0\}$ of stabilities for ${\bf d}$ and its associated real vector space ${\rm Stab}({\bf d})_\mathbb{R}={\rm Stab}({\bf d})\otimes_\mathbb{Z}\mathbb{R}$. For every non-zero proper ${\bf e}\leq{\bf d}$, we consider the hyperplane $W_{\bf e}=\{\Theta\, :\, \Theta({\bf e})=0\}\subset {\rm Stab}({\bf d})_\mathbb{R}$, called a wall in ${\rm Stab}({\bf d})_\mathbb{R}$. For a connected component $C_\mathbb{R}$ of the complement $${\rm Stab}^0({\bf d})_\mathbb{R}={\rm Stab}({\bf d})_\mathbb{R}\setminus\bigcup_{\bf e}W_e,$$
its closure $\overline{C_\mathbb{R}}$, resp.~the set of integral points $\overline{C_\mathbb{R}}\cap{\rm Stab}({\bf d})$ in it, is a (convex polyhedral) cone, called a chamber in ${\rm Stab}({\bf d})_\mathbb{R}$ and in ${\rm Stab}({\bf d})$, respectively.  Then ${\bf d}$ is $\Theta$-coprime if and only if $\Theta$ belongs to the interior of a chamber. If ${\bf d}$ is $\{{\bf d},\_\}$-coprime, we denote by $C_{\rm can}$ the chamber whose interior contains $\{{\bf d},\_\}$ and call it the canonical chamber.\\[1ex]
For a non-zero stability $\Theta\in{\rm Stab}({\bf d})$, we finally define
$$\gcd(\Theta)=\gcd(\Theta_i\, |\, i\in Q_0).$$

\section{Line bundles}\label{s2}

\subsection{Construction of tautological bundles on quiver moduli}

We review the construction of tautological bundles of \cite{King}. Assuming again that ${\bf d}$ is $\Theta$-coprime, and thus that ${\bf d}$ is indivisible, we can choose integers $a_i$ for $i\in Q_0$ such that $$\sum_ia_id_i=1,$$ and define $a\in(\mathbb{Z}Q_0)^*$ by $a({\bf i})=a_i$.\\[1ex]
We consider the trivial bundle $R_{\bf d}^{\Theta-{\rm st}}(Q)\times V_i\rightarrow R_{\bf d}^{\Theta-{\rm st}}(Q)$ and $G_{\bf d}$-linearize it by defining the action on $V_i$ by $$(g_j)_j\cdot v_i=\prod_{j\in Q_0}\det(g_j)^{-a_j}g_iv_i.$$
By definition of $a$, the scalar subgroup $\mathbb{C}^*$ then acts trivially, allowing this bundle to descend via the quotient map $R_{\bf d}^{\Theta-{\rm st}}(Q)\rightarrow M_{\bf d}^{\Theta-{\rm st}}(Q)$ to a bundle $\mathcal{V}_i$. By slight abuse of notation, we denote the above trivial bundle with its $G_{\bf d}$-linearization just by $V_i$; thus $V_i$ descends to $\mathcal{V}_i$. The bundles $\mathcal{V}_i$ are then tautological in the following sense:\\[1ex]
For any arrow $\alpha:i\rightarrow j$, there is a natural map $V_\alpha:V_i\rightarrow V_j$ given by

$$V_\alpha((f_\beta)_\beta,v_i))=((f_\beta)_\beta,f_\alpha(v_i)).$$

It descends to a map of vector bundles $\mathcal{V}_\alpha:\mathcal{V}_i\rightarrow\mathcal{V}_j$, thus defining a representation $\mathcal{V}$ of $Q$ in the category of vector bundles on $M_{\bf d}^{\Theta-{\rm st}}(Q)$, such that the induced quiver representation 
$$(((\mathcal{V}_i)_{[V]})_i,((\mathcal{V}_\alpha)_{[V]})_\alpha)$$
in the fibre over a point $[V]$ of $M_{\bf d}^{\Theta-{\rm st}}$ is isomorphic to $V$. This in fact makes $M_{\bf d}^{\Theta-{\rm st}}(Q)$ a fine moduli space, in the sense that any fiberwise $\Theta$-stable representation of $Q$ in vector bundles of rank vector ${\bf d}$ arises via pullback from the above tautological family.

\subsection{The Picard group of a fine quiver moduli space}

Assume as before that ${\bf d}$ is indivisible. The exact sequence of groups 
$$1\rightarrow \mathbb{C}^*\rightarrow G_{\bf d}\rightarrow PG_{\bf d}\rightarrow 1$$
induces an exact sequence
$$0\rightarrow X(PG_{\bf d})\rightarrow X(G_{\bf d})\rightarrow X(\mathbb{C}^*)$$
of character groups. The isomorphism 
$$(\mathbb{Z}Q_0)^*\rightarrow X(G_{\bf d})$$
mapping $\Theta$ to the character $\chi_\Theta$ as defined above yields a diagram of abelian groups
\enlargethispage{1em}
$$
	\begin{tikzcd}
		0 \arrow{r}{} &[-.5em] {\rm Stab}({\bf d}) \arrow{r}{} \arrow{d}{\simeq} & (\mathbb{Z}Q_0)^* \arrow{r}{\mathbf{d}} \arrow{d}{\simeq} & \mathbb{Z} \arrow{r}{} \arrow{d}{\simeq} &[-.5em] 0 \\[-.5em]
		0 \arrow{r}{} & X(PG_{\bf d}) \arrow{r}{} & X(G_{\bf d}) \arrow{r}{} & X(\mathbb{C}^*) \arrow{r}{} & 0
	\end{tikzcd}
$$
with right exact sequences since the evaluation at ${\bf d}$ admits a section induced by the element $a\in(\mathbb{Z}Q_0)^*$ chosen above.\\[1ex]
We have an isomorphism
$$X(PG_{\bf d})\simeq{\rm Pic}^{PG_{\bf d}}(R_{\bf d}(Q))$$
given by assigning to a character $\chi$ the trivial line bundle on $R_{\bf d}(Q)$ with $PG_{\bf d}$-linearization given by $\chi$, denoted $L(\chi)$.\\[1ex]
We now make the following additional assumption on ${\bf d}$:\\[1ex]
The dimension vector ${\bf d}$ is called {\it $\Theta$-amply stable} if the codimension of the unstable locus is at least two, that is
$${\rm codim}_{R_{\bf d}(Q)}(R_{\bf d}(Q)\setminus R_{\bf d}^{\Theta-{\rm st}}(Q))\geq 2.$$
A sufficient numerical condition for this property is stated in \cite[Proposition 5.1]{RS}: ${\bf d}$ is $\Theta$-amply stable if $\langle{\bf e},{\bf d}-{\bf e}\rangle\leq -2$ for every non-zero proper ${\bf e}\leq{\bf d}$ such that $\Theta({\bf e})\geq 0$.\\[1ex]
In this case, the natural restriction map
$${\rm Pic}^{PG_{\bf d}}(R_{\bf d}(Q))\rightarrow{\rm Pic}^{PG_{\bf d}}(R_{\bf d}^{\Theta-{\rm st}}(Q))$$
is an isomorphism. But since $R_{\bf d}^{\Theta-{\rm st}}(Q)$ is a $PG_{\bf d}$-principal bundle over $M_{\bf d}^{\Theta-{\rm st}}(Q)$, we have an isomorphism
$${\rm Pic}^{PG_{\bf d}}(R_{\bf d}^{\Theta-{\rm st}}(Q))\simeq{\rm Pic}(M_{\bf d}^{\Theta-{\rm st}}(Q)),$$
thus every line bundle $L(\chi)$ descends to a line bundle $\mathcal{L}(\chi)$ on $M_{\bf d}^{\Theta-{\rm st}}(Q)$.
Furthermore, since $M_{\bf d}^{\Theta-{\rm st}}(Q)$ is smooth, taking the first Chern class gives an isomorphism
$${\rm Pic}(M_{\bf d}^{\Theta-{\rm st}}(Q))\simeq A^1(M_{\bf d}^{\Theta-{\rm st}}(Q)).$$
We conclude:

\begin{proposition} If ${\bf d}$ is $\Theta$-coprime and $\Theta$-amply stable, we have a chain of isomorphisms of abelian groups
$$\begin{array}{ccccccc}{\rm Stab}({\bf d})&\stackrel{\simeq}{\rightarrow}&X(PG_{\bf d})&\stackrel{\simeq}{\rightarrow}&{\rm Pic}(M_{\bf d}^{\Theta-{\rm st}}(Q))&\stackrel{\simeq}{\rightarrow}&A^1(M_{\bf d}^{\Theta-{\rm st}}(Q))\\
\Theta'&\mapsto&\chi_{\Theta'}&\mapsto&\mathcal{L}({\chi_{\Theta'}})&\mapsto&c_1(\mathcal{L}({\chi_{\Theta'}})).\end{array}$$
\end{proposition}

We proceed with a stability condition $\Theta$ for which $\mathbf{d}$ is $\Theta$-coprime and $\Theta$-amply stable.
We want to identify the determinant bundles $\det(\mathcal{V}_i)$ of the tautological bundles $\mathcal{V}_i$ introduced above under this identification. By the definition of the $G_{\bf d}$-linearization of the bundle $\mathcal{V}_i$ on $R_{\bf d}^{\Theta-{\rm st}}(Q)$, its determinant bundle is linearized by the character
$$\chi((g_j)_j)=\prod_j(\det g_j)^{-a_jd_i}\cdot\det(g_i).$$
Denoting by $r:(\mathbb{Z}Q_0)^*\rightarrow{\rm Stab}({\bf d})$ the retraction induced by the section $a$, that is,
$$r(\Theta)=\Theta-\Theta({\bf d})\cdot a,$$
we thus have:

\begin{lemma} For all $i\in Q_0$, the following equality holds in ${\rm Pic}(M_{\bf d}^{\Theta-{\rm st}}(Q))$:
$$\det(\mathcal{V}_i)=\mathcal{L}({\chi_{-r(i)}}).$$
\end{lemma}

\subsection{The ample cone}

By slight abuse of notation, we denote the total space of a line bundle by the same symbol. For any $\chi\in X(PG_{\bf d})$, we then have a commutative square with the horizontal maps being $PG_{\bf d}$-principal bundles and the vertical maps being bundle projections
$$
	\begin{tikzcd}
		L(\chi) \arrow{d}{} \arrow{r}{} & \mathcal{L}(\chi) \arrow{d}{} \\[-.5em]
		R_{\bf d}^{\Theta-{\rm st}}(Q) \arrow{r}{} & M_{\bf d}^{\Theta-{\rm st}}(Q).
	\end{tikzcd}
$$
A global section of $\mathcal{L}(\chi)$ then corresponds uniquely to a $PG_{\bf d}$-equivariant section $R_{\bf d}^{\Theta-{\rm st}}(Q)\rightarrow L(\chi)$ of $L(\chi)$, which determines, and is determined by, a $\chi$-semi-invariant function on $R_{\bf d}^{\Theta-{\rm st}}(Q)$. This proves:
\begin{lemma} If ${\bf d}$ is $\Theta$-coprime and $\Theta$-amply stable, we have
$$H^0(M_{\bf d}^{\Theta-{\rm st}}(Q),\mathcal{L}(\chi))\simeq\mathbb{C}[R_{\bf d}(Q)]^{G_{\bf d},\chi}.$$
\end{lemma}

We show the following characterization of the ample cone of the moduli space $M_{\bf d}^{\Theta-{\rm st}}(Q)$:

\begin{proposition} If ${\bf d}$ is $\Theta$-coprime and $\Theta$-amply stable, and if $\Theta'\in{\rm Stab}({\bf d})$ belongs to the interior of the same chamber as $\Theta$, the line bundle $\mathcal{L}(\chi_{\Theta'})$ on $M_{\bf d}^{\Theta-{\rm st}}(Q)$ is ample. In particular, $\mathcal{L}(\chi_\Theta)$ is ample. 
\end{proposition}

\proof The line bundle $\mathcal{L}(\chi_{\Theta'})$ being ample means that, for large enough $N$, generators $f_0,\ldots,f_s$ of $$H^0(M_{\bf d}^{\Theta-{\rm st}}(Q),\mathcal{L}(\chi_{\Theta'})^{\otimes N})\simeq\mathbb{C}[R_{\bf d}(Q)]^{G_{\bf d},N\chi_{\Theta'}}$$
descend from $R_{\bf d}^{\Theta-{\rm st}}(Q)$ via the quotient map to a closed embedding
$$(f_0:\ldots:f_s):M_{\bf d}^{\Theta-{\rm st}}(Q)\rightarrow\mathbb{P}^s.$$
But by the construction of quiver moduli, we know that, assuming $N$ to be large enough to generate $\mathbb{C}[R_{\bf d}(Q)]^{G_{\bf d}}_{\chi_{\Theta'}}$, the image of this map is precisely $M_{\bf d}^{\Theta'-{\rm st}}(Q)$. Ampleness of $\mathcal{L}(\chi_{\Theta'})$ is thus
equivalent to such functions inducing a closed immersion
$$M_{\bf d}^{\Theta-{\rm st}}(Q)\rightarrow M_{\bf d}^{\Theta'-{\rm st}}(Q),$$
and thus an isomorphism since both spaces are irreducible and projective. This is certainly the case if the open sets $R_{\bf d}^{\Theta-{\rm st}}(Q)$ and $R_{\bf d}^{\Theta'-{\rm st}}(Q)$ coincide, which in turn is fulfilled if $\Theta$ and $\Theta'$ belong to the interior of the same chamber in ${\rm Stab}({\bf d})$.\\[1ex]
We have thus identified the interior of the chamber of a stability $\Theta$, for which ${\bf d}$ is $\Theta$-coprime and $\Theta$-amply stable, as part of the ample cone of $M_{\bf d}^{\Theta-{\rm st}}(Q)$.

\section{Fano quiver moduli}\label{s3}

\subsection{Computation of the (anti-)canonical class}


We would like to express the class of the determinant bundle $\det\mathcal{T}$ in ${\rm Pic}(M_{\bf d}^{\Theta-{\rm st}}(Q))$, resp.~in $A^1(M_{\bf d}^{\Theta-{\rm st}}(Q))$, where $\mathcal{T}$ is the tangent bundle. To do this, note that the fiber of $\mathcal{T}$ in a point $[V] \in M_{\bf d}^{\Theta-{\rm st}}(Q)$ is
$$
\mathcal{T}_{[V]} \simeq {\rm Ext}^1_Q(V,V).
$$
Let $[V]$ be represented by $V = (f_\alpha)_\alpha \in R_{\bf d}(Q)$. There is an exact sequence
$$
	0 \to \End_Q(V,V) \to \bigoplus_{i \in Q_0} \End(V_i) \xrightarrow[]{\phi} \bigoplus_{\alpha: i \to j} \Hom(V_i,V_j) \to \Ext_Q^1(V,V) \to 0
$$
where $\phi((x_i)_i) = (x_jf_\alpha - f_\alpha x_i)_\alpha$. This sequence is a consequence of the standard projective resolution of $V$. Note that the map $\phi$ is precisely the derivative of the map $G_{\bf d} \to R_{\bf d}^{\Theta-{\rm st}}(Q)$ which sends $g \mapsto gV$. The endomorphism ring of a stable representation reduces to the scalars. By globalizing the above exact sequence, we obtain an exact sequence of vector bundles
$$
0\rightarrow \mathcal{O} \rightarrow \bigoplus_{i\in Q_0} \mathcal{V}_i^* \otimes \mathcal{V}_i \rightarrow \bigoplus_{\alpha:i\rightarrow j} \mathcal{V}_i^* \otimes \mathcal{V}_j \rightarrow \mathcal{T} \rightarrow 0.
$$

To calculate the anticanonical class, we perform a Chern class calculation in the Chow group $A^1(M_{\bf d}^{\Theta-{\rm st}}(Q))$, made possible by the above isomorphism to the Picard group, using the following properties of first Chern classes:
\begin{itemize}
\item $c_1$ is additive on short exact sequences
\item $c_1(\det\mathcal{E})=c_1(\mathcal{E})$,
\item $c_1(\mathcal{O})=0$,
\item $c_1(\mathcal{E}^*)=-c_1(\mathcal{E})$,
\item $c_1(\mathcal{E}\otimes\mathcal{F})={\rm rk}(\mathcal{F})\cdot c_1(\mathcal{E})+{\rm rk}(\mathcal{E})\cdot c_1(\mathcal{F})$.
\end{itemize}
We thus find
\begin{align*}
c_1(\bigoplus_{i\in Q_0} \mathcal{V}_i^* \otimes \mathcal{V}_i ) &= \sum_i(-d_i c_1(\mathcal{V}_i)+d_i c_1(\mathcal{V}_i))=0,\\
c_1(\bigoplus_{\alpha:i\rightarrow j} \mathcal{V}_i^* \otimes \mathcal{V}_j) &= \sum_{\alpha:i\rightarrow j}(-d_j c_1(\mathcal{V}_i)+d_i c_1(\mathcal{V}_j)),
\end{align*}
and thus, using the above exact sequence,
\begin{align*}
c_1(\mathcal{T}) &= c_1(\bigoplus_{\alpha:i\rightarrow j}\mathcal{V}_i^* \otimes \mathcal{V}_j) \\
&= \sum_{\alpha:i\rightarrow j}(-d_j c_1(\mathcal{V}_i)+d_i c_1(\mathcal{V}_j)) \\
&= \sum_i(\langle {\bf i},{\bf d}\rangle-\langle {\bf d},{\bf i}\rangle)c_1(\mathcal{V}_i).
\end{align*}
Since 
$$c_1(\mathcal{V}_i)=c_1(\det(\mathcal{V}_i))=c_1(\mathcal{L}(\chi_{-r(i)}),$$
we thus find $\det(\mathcal{T})=\mathcal{L}(\chi_{\Theta'})$ for
$$\Theta'=\sum_i(\langle{\bf i},{\bf d}\rangle-\langle{\bf d},{\bf i}\rangle)(-r(i))=-r(\langle\_,{\bf d}\rangle-\langle{\bf d},\_\rangle)=r(\{{\bf d},\_\})=\{{\bf d},\_\}$$
since $r$ acts as the identity on linear forms in ${\rm Stab}({\bf d})$. We have thus proved:
\begin{proposition} 
In ${\rm Pic}(M_{\bf d}^{\Theta-{\rm st}})(Q)$, we have $$\det(\mathcal{T})=\mathcal{L}(\chi_{\{{\bf d},\_\}}).$$
\end{proposition}

Applying our above criterion for ampleness, we can thus derive our main result, summarizing our findings:

\begin{theorem}\label{mainthm} If ${\bf d}$ is $\{{\bf d},\_\}$-coprime and $\{{\bf d},\_\}$-amply stable, then the moduli space $M_{\bf d}^{\{{\bf d},\_\}-{\rm st}}(Q)$ is a smooth irreducible projective Fano variety of dimension $1-\langle{\bf d},{\bf d}\rangle$, Picard rank $|Q_0|-1$ and index $\gcd(\{{\bf d}\,\_\})$, which is rational and has only  algebraic cohomology.
\end{theorem}

\section{Classes of examples}\label{s4}

\subsection{Subspace quivers -- moduli of point configurations in projective space}

We consider the quiver $S_m$ with vertices $i_1,\ldots,i_m,j$ and arrows $i_k\rightarrow j$ for $k=1,\ldots,m$, called the $m$-subspace quiver. For a fixed $d\geq 2$, we consider the dimension vector $${\bf d}=\sum_k{\bf i}_k+d{\bf j}.$$ We have $$\Theta=\{{\bf d},\_\}=d\sum_{k}i_k-mj,$$ thus ${\bf d}$ is $\Theta$-coprime if $m$ and $d$ are coprime. In this case, the moduli space $M_{\bf d}^{\Theta-st}(S_m)$ equals the moduli space of stable ordered configurations of $m$ points in $\mathbb{P}^{m-1}$ modulo the natural ${\rm PGL}_m(\mathbb{C})$-action. Here $(p_1,\ldots,p_m)$ is called stable if the linear subspace spanned by any $k$ of these points has dimension strictly bigger than $mk/d-1$. This moduli space is non-empty and not a single point if and only if $m-1>d$, in which case it is of dimension $$(d-1)(m-d-1).$$ We can easily verify that ${\bf d}$ is $\{{\bf d},\_\}$-amply stable with the numerical sufficient criterion: suppose that ${\bf e}\leq{\bf d}$, thus ${\bf e}=\sum_k{\bf i}_k+e{\bf j}$ for a subset $K\subset\{1,\ldots,m\}$ and $e\leq d$. Then $\Theta({\bf e})>0$ if and only if $|K|>me/d$. We calculate
$$\langle{\bf e},{\bf d}-{\bf e}\rangle_{S_m}=(e-|K|)(d-e)<-\frac{m-d}{d}e(d-e)\leq 0.$$
Thus, assuming that $\langle{\bf e},{\bf d}-{\bf e}\rangle\geq -1$, we find $$e=d-1\mbox{ and }|K|=e+1=d.$$ By assumption on $K$, this yields $d>m(d-1)/d$ or, equivalently, $$d^2>m(d-1)\geq d^2-1.$$ We thus find $m=d+1$, contradicting the assumption.

\begin{corollary} For given coprime $m$ and $d$ such that $m-1>d\geq 2$, the moduli space $(\mathbb{P}^{m-1})^d_{\rm st}/{\rm PGL}_d(\mathbb{C})$ is Fano of dimension $(d-1)(m-d-1)$, with Picard rank $m$ and index one.
\end{corollary}

To examine these Fano varieties in low dimension, we can additionally assume $m>2d$ using the duality \cite[Theorem 3.1]{RW}. For $d=2$ and $m=5$ we thus find a del Pezzo surface of Picard rank $5$, which is therefore isomorphic to the blowup of the projective plane in four general points. We will see below how all del Pezzo surfaces of lower Picard rank (which are already toric) can be realized as quiver moduli. For $d=2$ and $m=7$ we find a four-dimensional Fano variety with Picard rank $7$ and fourth Betti number $b_4=22$. For $d=3$ and $m=7$ we find a six-dimensional Fano variety with Picard rank $7$ and $b_4=29$, $b_6=64$, using the Betti number formulas of \cite{RW}.\\[1ex]
An interesting (non-)example arises here, too, which shows that the assumption of our main result cannot be relaxed much: for the moduli space of six points in the projective line, the dimension vector ${\bf d}$ is not $\Theta$-coprime for the canonical stability $\Theta$. The semistable moduli space is isomorphic to the Segre cubic by \cite{HMSV}, which is indeed a singular projective Fano threefold with ten isolated singularities. Slightly deforming the stability $\Theta$ to a new stability $\Theta^+$ as in \cite{FR} yields a small desingularization of the Segre cubic, which cannot be Fano: a Fano threefold of Picard rank six is automatically isomorphic to the product of the projective line and the degree five del Pezzo surface, thus contains a $2$-nilpotent element in second cohomology. But in \cite{FR} it is shown that this does not hold for this desingularization, disproving the Fano property.

\subsection{Generalized Kronecker quivers -- Kronecker moduli}

Let $K_m$ be the quiver with vertices $i$ and $j$ and $m\geq 1$ arrows from $i$ to $j$. Let ${\bf d}=d{\bf i}+e{\bf j}$ be a dimension vector, and let $$\Theta=\{{\bf d},\_\}=mei-mdj$$ be the canonical stability. As in the previous example, $\Theta$-coprimality of ${\bf d}$ is equivalent to $d$ and $e$ being coprime. If $m=1$ or $m=2$, the only non-empty moduli spaces $K_{(d,e)}^{(m)}=M_{\bf d}^{\Theta-{\rm st}}(Q)$ are single points or a projective line (for $m=2$ and $d=1=e$), thus we assume $m\geq 3$. In this case, if non-empty, the moduli space is of dimension $$mde-d^2-e^2+1.$$ The $\Theta$-ample stability of ${\bf d}$ is proved in \cite[Proposition 6.2]{RS}. We thus find:

\begin{corollary} For given $m\geq 3$ and coprime $d$ and $e$, the Kronecker moduli space $K_{(d,e)}^{(m)}$, if non-empty, is a Fano variety of dimension $mde-d^2-e^2+1$, Picard rank one, and index $m$.
\end{corollary}

A theorem of Kobayashi and Ochiai \cite{KO} asserts that a smooth projective Fano variety $X$ of dimension $n$ and index $q$ satisfies $q \leq n+1$ and equality holds if and only if $X \simeq \mathbb{P}^n$. This led Mukai \cite{Muk} to a conjecture on the relation between the rank, the index and the dimension of a Fano variety. This conjecture was later refined by Bonavero, Casagrande, Debarre and Druel \cite{BCDD}. It states the following. Let $X$ be a smooth projective Fano variety of dimension $n$, let $r$ denote the rank of the Picard group of $X$, and let $p$ be the minimal positive integer such that $-K_X \cdot C = p$ for a rational curve $C \subseteq X$, called the pseudo-index of $X$. It is then conjectured that
$$
r\cdot (p-1) \leq n,
$$
and equality is satisfied if and only if $X$ is isomorphic to $(\mathbb{P}^{p-1})^r$. In the original Mukai conjecture, the pseduo-index is replaced by the index. Note that $q\leq p$, and by the Kobayashi theorem together with $p\leq n+1$, proved by Mori \cite{Mori}, one obtains $q\leq p\leq n+1$. It would be interesting to know if the Mukai conjecture holds for Fano quiver moduli (see the following section for a first example). For Kronecker moduli, the rank of the Picard group is one, so the conjecture is true by the Kobayashi--Ochiai theorem. We illustrate how to confirm this statement entirely numerically in this case.

\begin{proposition} If non-empty, the moduli space $K_{(d,e)}^{(m)}$ has at least dimension $m-1$, and is isomorphic to $\mathbb{P}^{m-1}$ if equality holds.
\end{proposition}

\proof So assume that $K_{(d,e)}^{(m)}$ is non-empty. We will prove that $$mde-d^2-e^2+1\geq m-1,$$ and that equality holds if and only if $d=1=e$. By standard reductions (see the proof of \cite[Proposition 6.2]{RS}) we can assume that $d\leq e\leq md/2$. Writing $e=\alpha d$ for $\alpha\in\mathbb{Q}$ with $1\leq\alpha\leq m/2$, we thus have to prove that $$(m\alpha-1-\alpha^2)d^2\geq m-2,$$ and that equality holds if and only if $d=1=e$. The quadratic function $$\alpha\mapsto m\alpha-1-\alpha^2$$ being strictly monotonously increasing in the interval $[1,m/2[$, we have $$m\alpha-1-\alpha^2\geq m-2,$$ with equality only for $\alpha=1$, and thus $$(m\alpha-1-\alpha^2)d^2\geq (m-2)d^2\geq m-2,$$ with equality only if $d=1$ and $\alpha=1$, which implies $d=1=e$ by coprimality of $d$ and $e$.

\subsection{A class of Fano quiver moduli of arbitrary rank and index}

For fixed $m,k\geq 1$, we consider the quiver $S^{(k)}_m$ with vertices $i_1,\ldots,i_m,j$ and $k$ arrows from each vertex $i_l$ to $j$, for $l=1,\ldots,m$ (a thickened version of the $m$-subspace quiver). For a fixed $d\geq 1$, we consider the dimension vector $${\bf d}=\sum_k{\bf i}_k+d{\bf j}.$$ 
We have $$\Theta=\{{\bf d},\_\}=k(d\sum_{l}i_l-mj),$$ thus ${\bf d}$ is again $\Theta$-coprime if $m$ and $d$ are coprime, which we assume from now on.

A representation of $S^{(k)}_m$ is given by a tuple $(v_{i,j}\, |\, i=1,\ldots,m,\, j=1,\ldots, k)$ of vectors in $\mathbb{C}^d$, which are considered up to the diagonal action of ${\rm GL}_d(\mathbb{C)}$, and up to the action of an $m$-torus $T$ acting by
$$(t_i)_i\cdot(v_{i,j})_{i,j}=(t_iv_{i,j})_{i,j}.$$

Such a tuple is stable if the $v_{i,j}$ span $\mathbb{C}^d$ and, for every nonempty proper subset $I$ of $\{1,\ldots,m\}$, the span $U_I$ of the $v_{i,j}$ for $i\in I$ has dimension strictly bigger than $d|I|/m$. Assuming $d\leq km$, we can choose the vectors $v_{i,j}$ sufficiently independent such that, for all nonempty proper $I$, the space $U_I$ has dimension $k|I|$ if $k|I|\leq d$, and equals $\mathbb{C}^d$ otherwise. In both cases, the estimate for stability is trivially fulfilled, proving that there exists a stable point if $d\leq km$. 

Arranging the vectors $v_{i,j}$ into a $d\times km$-matrix and noting that the quotient of the open set of highest rank matrices by the ${\rm GL}_d(\mathbb{C})$-action is isomorphic to the Grassmannian ${\rm Gr}_d^{km}$, we see that the moduli space $M_{\bf d}^{\Theta-st}(S^{(k)}_m)$ admits an interpretation as the quotient of an open set of stable points in ${\rm Gr}_d(\mathbb{C}^k\otimes\mathbb{C}^m)$ by the action of the torus $T$ induced from its natural action on $\mathbb{C}^m$. 

As special cases, we find the moduli spaces of point configurations from above for $k=1$, Grassmannians ${\rm Gr}_d^k$ for $m=1$, and powers of projective spaces $(\mathbb{P}^{k-1})^m$ for $d=1$ or $d=mk-1$.

In general, the moduli space is of dimension
$$1-m-d^2+kmd=(km-1-d)(d-1)+(k-1)m,$$
thus it is not empty and not a single point if $d<km$, and additionally $2\leq d\leq m-1$ in the case $k=1$ treated before; we assume these estimates for $d$ from now on.

We can then verify $\Theta$-ample stability: suppose given a dimension vector $\mathbf{e}\leq\mathbf{d}$ such that $\Theta(\mathbf{e})>0$ and $\langle\mathbf{e},\mathbf{d}-\mathbf{e}\rangle\geq -1$. Then $\mathbf{e}=\sum_{k\in I}{\bf i}_k+e{\bf j}$ for a subset $I$ of $\{1,\ldots,m\}$ and an $e\leq d$. Denoting $n=|I|$, the assumptions on $\mathbf{e}$ translate to $dn>em$ and $(kn-e)(d-e)\leq 1$. The first inequality $dn>em$ easily implies that both factors in the second inequality are positive, and thus $kn-1=e=d-1$. Then $dn>em$ reads $d\cdot d/k>(d-1)m$, and thus $d^2>(d-1)mk$, which is equivalent to
$$(d-1)(mk-1-d)\leq 0.$$
Both factors being nonnegative, we arrive at $d=1$ or $d=mk-1$. But we also have $k|d$, and thus $k=1$, which is the trivial case we excluded above. We thus find:

\begin{proposition}\label{vmukai} For $m,k,d\geq 1$ such that $d$ and $m$ are coprime, $d\leq mk-1$, and $d\not=1,m-1$ in case $k=1$, the moduli space ${\rm Gr}_d(\mathbb{C}^k\otimes\mathbb{C}^m)/T$ is a Fano variety of dimension $$(km-1-d)(d-1)+(k-1)m,$$
Picard rank $m$ and index $k$, and it verifies the Mukai conjecture.
\end{proposition}

\proof It remains to verify the Mukai conjecture: our dimension formula indeed yields the estimate
$$m(k-1)\leq(km-1-d)(d-1)+(k-1)m,$$
and equality holds if $d=1$ or $d=km-1$, in which case the moduli space is a power of projective space.

 \subsection{Toric quiver moduli}

On the variety of representations $R_{\bf d}(Q)$, there is an action of a torus $T_1$ of rank equal to the number of arrows in $Q$ by scaling the linear maps,
$$(t_\alpha)_\alpha(f_\alpha)_\alpha=(t_\alpha f_\alpha)_\alpha,$$
which commutes with the base change action of $G_{\bf d}$, thus inducing an action on $M_{\bf d}^{\Theta-{\rm st}}(Q)$. If ${\bf d}={\bf 1}:=\sum_{i\in Q_0}{\bf i}$, this endows $M_{\bf 1}^{\Theta-{\rm st}}(Q)$ with the structure of a toric variety. We make our conditions for the moduli space being Fano explicit in this case.\\[1ex]
Since $Q$ is acyclic, we can assume that $Q_0=\{i_1,\ldots,i_n\}$ is indexed such that existence of an arrow $\alpha:i_k\rightarrow i_l$ implies $k<l$. We denote by $a_{k,l}$ the number of arrows from $k$ to $l$. For subsets $K,L\subset [n]=\{1,\ldots,n\}$, we abbreviate
$$a_{K,L}=\sum_{k\in K}\sum_{l\in L}a_{k,l}.$$ The canonical stability $\Theta=\{{\bf 1},\_\}$ is given by
$$\Theta_i=a_{i[n]}-a_{[n]i}.$$
$\Theta$-coprimality then easily translates into the condition $a_{K[n]}\not=a_{[n]K}$ for all proper non-empty subsets $K$, which simplifies to
$$a_{K\overline{K}}\not=a_{\overline{K}K}$$
for all proper non-empty $K$, where $\overline{K}=[n]\setminus K$. Finally, the numerical condition for ample stability reads
$$\max(a_{K\overline{K}},a_{\overline{K}K})\geq 2.$$
We thus find:

\begin{proposition} Given nonnegative integers $a_{kl}$ for $1\leq k<l\leq n$ such that
$$a_{K\overline{K}}\not=a_{\overline{K}K}\mbox{ and }\max(a_{K\overline{K}},a_{\overline{K}K})\geq 2$$
for all proper nonempty subsets $K\subset[n]$, the moduli space $M_{\bf 1}^{\{ {\bf 1},\_\}-{\rm st}}(Q)$ for the quiver with vertices $i_1,\ldots,i_n$ and $a_{kl}$ arrows from vertex $i_k$ to $i_l$ is a smooth projective toric Fano variety of dimension $a_{[n][n]}-n+1$, Picard rank $n-1$ and index $\gcd(a_{k[n]}-a_{[n]k}\, |\, k=1,\ldots,n)$.
\end{proposition}

With this result, we can find the del Pezzo surfaces $\mathbb{P}^1\times\mathbb{P}^1$, ${\rm Bl}_1\mathbb{P}^2$, ${\rm Bl}_2\mathbb{P}^2$ and ${\rm Bl}_3\mathbb{P}^2$ by the respective toric quiver moduli
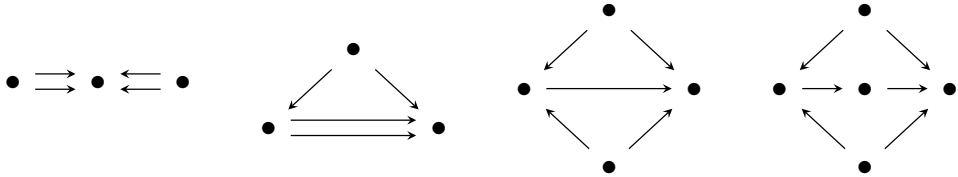
\tikzcdset{arrow style=tikz, diagrams={>=stealth}}
\tikzset{commutative diagrams/row sep/normal=3ex}
\tikzset{commutative diagrams/column sep/normal=3ex}
\begin{center}
	$\begin{tikzcd}
		\bullet \arrow[shift right]{r}{} \arrow[shift left]{r}{} & \bullet & \arrow[shift right]{l}{} \arrow[shift left]{l}{}\bullet
	\end{tikzcd}$
\quad
	$\begin{tikzcd}
		& \bullet \arrow{dr}{} \arrow{dl}{} & \\
		\bullet \arrow[shift right]{rr}{} \arrow[shift left]{rr}{}& & \bullet
	\end{tikzcd}$
\quad
	$\begin{tikzcd}
		& \bullet \arrow{dr}{} \arrow{dl}{} & \\
		\bullet \arrow{rr}{} & & \bullet \\
		& \bullet \arrow{ur} \arrow{ul}{} &
	\end{tikzcd}$
\quad
	$\begin{tikzcd}
		& \bullet \arrow{dr}{} \arrow{dl}{} & \\
		\bullet \arrow{r}{} & \bullet \arrow{r}{} & \bullet \\
		& \bullet \arrow{ur}{} \arrow{ul}{} &
	\end{tikzcd}$
\end{center}

We can also realize the toric Fano threefolds of Picard rank two $\mathbb{P}^1\times\mathbb{P}^2$, ${\rm Bl}_p\mathbb{P}^3$ (blowup in a point), ${\rm Bl}_l\mathbb{P}^3$ (blowup in a line), respectively, by
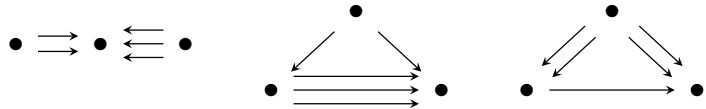
\begin{center}
	$\begin{tikzcd}
		\bullet \arrow[shift right]{r}{} \arrow[shift left]{r}{} & \bullet & \arrow[shift right=1ex]{l}{} \arrow{l}{} \arrow[shift left=1ex]{l}{}\bullet
	\end{tikzcd}$
\quad
	$\begin{tikzcd}
		& \bullet \arrow{dr}{} \arrow{dl}{} & \\
		\bullet \arrow[shift right=1ex]{rr}{} \arrow{rr}{} \arrow[shift left=1ex]{rr}{} & & \bullet
	\end{tikzcd}$
\quad
	$\begin{tikzcd}
		& \bullet \arrow[shift left]{dr}{} \arrow[shift right]{dr}{} \arrow[shift left]{dl}{} \arrow[shift right]{dl}{}& \\
		\bullet \arrow{rr}{} & & \bullet
	\end{tikzcd}$
\end{center}

The two blowups can be distinguished by considering the index; we note that a fourth class of such Fano varieties, namely $\mathbb{P}(\mathcal{O}_{\mathbb{P}^2(2)}\oplus\mathcal{O}_{\mathbb{P}^2})$, cannot be realized using quiver moduli.


\begin{thebibliography}{99}

\bibitem{AH} K. Altmann, L. Hille, Strong exceptional sequences provided by quivers. Algebras Represent. Theory {\bf 2} (1999), 1--17.

\bibitem{BCDD} L. Bonavero, C. Casagrande, O. Debarre, S. Druel, Sur une conjecture de Mukai. Comment. Math. Helv. {\bf 78} (2003), no. 3, 601--626.
\bibitem{Fanog} P. Belmans, Fanography. \texttt{https://fanography.info}
\bibitem{Ca06} C. Casagrande, The number of vertices of a Fano polytope. Ann.
Inst. Fourier (Grenoble) {\bf 56} (2006), no. 1, 121-–130.


\bibitem{Fei} J. Fei, Moduli of Representations I. Projections from Quivers. \texttt{arXiv:1011.6106}

\bibitem{F1} H. Franzen, Chow rings of fine quiver moduli are tautologically presented. Math. Z. {\bf 279} (2015), no. 3--4, 1197-–1223.
\bibitem{F2} H. Franzen, 
Torus-Equivariant Chow Rings of Quiver Moduli. \texttt{arXiv:1911.03288}


\bibitem{FR} H. Franzen, M. Reineke, Cohomology rings of moduli of point configurations on the projective line. Proc. Amer. Math. Soc. {\bf 146} (2018), no. 6, 2327–-2341. 

\bibitem{GH} G. Gagliardi, J. Hofscheier, The generalized Mukai
conjecture for symmetric varieties. Trans. Amer. Math. Soc. {\bf 369} (2017),
2615--2649.

\bibitem{Hal} M. Halic, Strong exceptions sequences of vector bundles on certain Fano varieties. \texttt{arXiv:0906.3466}
\bibitem{Hil} L. Hille, Aktionen algebraischer Gruppen, geometrische Quotienten und K\"ocher. Habilitations\-schrift, Hamburg 2002.
\bibitem{HMSV} B. Howard, J. Millson, A. Snowden, R. Vakil,
A description of the outer automorphism of S6, and the invariants of six points in projective space. J. Combin. Theory Ser. A {\bf 115} (2008), no. 7, 1296-–1303. 
\bibitem{King} A. King,
Moduli of representations of finite-dimensional algebras.
Quart. J. Math. Oxford Ser. (2) {\bf 45} (1994), no. 180, 515-–530. 
\bibitem{KW} A. King, C. Walter,
On Chow rings of fine moduli spaces of modules.
J. Reine Angew. Math. {\bf 461} (1995), 179-–187. 
\bibitem{KO} S. Kobayashi, T. Ochiai, Characterization of complex projective spaces and hyperquadrics. J. Math. Kyoto Univ. {\bf 13} (1973), no. 1, 31--47. 
\bibitem{Mori} S. Mori, Projective manifolds with ample tangent bundles. Ann. of Math.
{\bf 110} (1979), 593--606.
\bibitem{Muk} S. Mukai, Problems on characterization of the complex projective space. Birational Geometry of Algebraic Varieties, Open Problems, Katata, 1988, 23rd Symposium, Taniguchi Foundation (1988), 57--60. 
\bibitem{Pas10} B. Pasquier, The pseudo-index of horospherical Fano varieties.
Internat. J. Math. {\bf 21} (2010), no. 9, 1147-–1156.

\bibitem{RHN} M. Reineke, The Harder-Narasimhan system in quantum groups and cohomology of quiver moduli.
Invent. Math. {\bf 152} (2003), no. 2, 349-–368. 
\bibitem{RS} M. Reineke, S. Schr\"oer,
Brauer groups for quiver moduli.
Algebr. Geom. {\bf 4} (2017), no. 4, 452–-471. 
\bibitem{RW} M. Reineke, T. Weist, Moduli spaces of point configurations and plane curve counts.
International Mathematical Research Notices, to appear.
\bibitem{Scho} A. Schofield, General representations of quivers. Proc. London Math. Soc. (3) {\bf 65} (1992), no. 1, 46-–64.
\bibitem{SchB} A. Schofield,
Birational classification of moduli spaces of representations of quivers.
 Indag. Math. (N.S.) {\bf 12} (2001), no. 3, 407-–432. 

\end{thebibliography}
\end{document}